\documentclass[reqno]{amsart}
\usepackage{amsmath, amssymb, amsthm, geometry, enumerate, graphicx, amsfonts, hyperref, mathrsfs}
\hypersetup{colorlinks=true,linkcolor=red, anchorcolor=blue, citecolor=blue, urlcolor=red, filecolor=magenta, pdftoolbar=true}
\theoremstyle{plain}
\newtheorem{thm}{\it Theorem}[section]

\newtheorem{cor}[thm]{\it Corollary}

\theoremstyle{remark}
\newtheorem{defn}[thm]{Def{}inition}
\newtheorem{rem}[thm]{Remark}

\numberwithin{equation}{section}
\usepackage{enumerate}
\allowdisplaybreaks
\begin{document}

\title []{On the paper ``Characterizations of weaving for $g$-frames by induced sequences''}

\author[Deepshikha]{ Deepshikha }
\address{{\bf{Deepshikha}}, Department of Mathematics, Shyampur Siddheswari  Mahavidyalaya \\
University of Calcutta, Kolkata, West Bengal-711312, India.}
\email{dpmmehra@gmail.com}

\author[Lalit   Kumar Vashisht]{ Lalit  Kumar  Vashisht}
\address{{\bf{Lalit  Kumar  Vashisht}}, Department of Mathematics,
University of Delhi, Delhi-110007, India.}
\email{lalitkvashisht@gmail.com}

\author[Geetika Verma]{Geetika Verma}
\address{{\bf{Geetika Verma}}, STEM College, School of Science\\
RMIT University, Melbourne Campus, Australia.}
\email{geetika.verma@rmit.edu.au}

\begin{abstract}
A counter-example by Xiangchun Xiao, Guoping Zhao and Guorong Zhou in [\textcolor[rgb]{0.00,0.07,1.00}{J. Pseudo-Differ. Oper. Appl.} (2021) 12:60] is incorrect. Further, there is no new characterization of weaving for $g$-frames by Xiangchun Xiao et al. in \cite{XZZ}. Xiangchun Xiao et al. in  \cite{XZZ} gave a wrong interpretation    about a  characterization of weaving $g$-frames in separable Hilbert spaces which already proved by Deepshikha, Vashisht and Verma  in \cite{DVV}.
\end{abstract}

\renewcommand{\thefootnote}{}


\maketitle
\baselineskip15pt
\section{Generalized Weaving Frames in Hilbert Spaces}
In \cite{DVV}, Deepshikha, Vashisht and Verma studied $\Theta$-$g$-woven frames in separable Hilbert spaces. They consider the following definition.
\begin{defn}\label{def1}\cite{DVV}
Two $\Theta$-frames $\{\Lambda_j^*f_{jk}\}_{j\in\mathbb{N}, k\in I_j}$  and $\{\Omega_j^*g_{jk}\}_{j\in\mathbb{N}, k\in Q_j}$ of  a separable Hilbert space $\mathcal{H}$ are said to be woven if there exist universal positive constants $A$ and $B$ such that for any $\sigma\subseteq\mathbb{N}$, the family $\{\Lambda_j^*f_{jk}\}_{j\in\sigma, k\in I_j}\cup\{\Omega_j^*g_{jk}\}_{j\in\mathbb{N}\setminus\sigma, k\in Q_j}$ is a $\Theta$-frame for $\mathcal{H}$ with lower frame bound $A$ and upper frame bound $B$.
\end{defn}
 Deepshikha et al. proved the following  characterization for $\Theta$-$g$-woven frames in terms of woven $\Theta$-frames in  \cite{DVV}.
\begin{thm}\label{thm1}\cite{DVV}
Suppose that $\Lambda \equiv\{\Lambda_{j}\}_{j=1}^\infty$  and  $\Omega \equiv\{\Omega_{j}\}_{j=1}^\infty$ are $\Theta$-$g$-frames for $\mathcal{H}$ with respect to $\{ \mathcal{H}_j\}_{j=1}^{\infty}$ and $\{ \mathcal{W}_j\}_{j=1}^{\infty}$, respectively. Let $\{f_{jk}\}_{k\in I_j \subset \mathbb{N}}$ and $\{g_{jk}\}_{k\in Q_j \subset \mathbb{N}}$ be frames for $\mathcal{H}_j$ and $\mathcal{W}_j$, respectively $(j \in \mathbb{N})$ with frame bounds $\alpha$, $\beta$ and $\alpha'$, $\beta'$, respectively.
Then the following conditions are equivalent.
\begin{enumerate}[$(i)$]
 \item $\Lambda$ and $\Omega$  are $\Theta$-$g$-woven.
 \item  $\{\Lambda_j^*f_{jk}\}_{j\in\mathbb{N}, k\in I_j}$ and $\{\Omega_j^*g_{jk}\}_{j\in\mathbb{N}, k\in Q_j}$ are woven $\Theta$-frames for $\mathcal{H}$.
 \end{enumerate}
\end{thm}
\begin{rem}In  the above theorem, we can see that the family  $\{\Lambda_j^*f_{jk}\}_{j\in\mathbb{N}, k\in I_j}$ has two indexing sets $\mathbb{N}$ and $I_j$, and the family $\{\Omega_j^*g_{jk}\}_{j\in\mathbb{N}, k\in Q_j}$ has two indexing sets $\mathbb{N}$ and $Q_j$. In general, the indexing sets $I_j$ and $Q_j$ are not same. So we can only take the partition of $\mathbb{N}$, the common indexing set. Thus, we consider Definition \ref{def1} (see, \cite{DVV} also) for the weaving of $\Theta$-frames $\{\Lambda_j^*f_{jk}\}_{j\in\mathbb{N}, k\in I_j}$  and $\{\Omega_j^*g_{jk}\}_{j\in\mathbb{N}, k\in Q_j}$  in Theorem \ref{thm1}.
\end{rem}
Xiangchun Xiao, Guoping Zhao and Guorong Zhou in  \cite{XZZ} consider the following definition for the weaving of $\Theta$-frames (or frames) $\{\Lambda_j^*f_{jk}\}_{j\in\mathbb{N}, k\in K_j}$  and $\{\Omega_j^*g_{jk}\}_{j\in\mathbb{N}, k\in K_j}$.
 \begin{defn}\label{def3}\cite{XZZ}
 Two $\Theta$-frames $\{\Lambda_j^*f_{jk}\}_{j\in\mathbb{N}, k\in K_j}$  and $\{\Omega_j^*g_{jk}\}_{j\in\mathbb{N}, k\in K_j}$ are said to be woven $\Theta$-frames for $\mathcal{H}$ if there exist universal positive constants $A$ and $B$ such that for any $\tau_0\subseteq\mathbb{N}$, and $\sigma_j\subseteq K_j$, $j\in\tau_0$, the family $\{\Lambda_j^*f_{jk}\}_{j\in\tau_0, k\in \sigma_j}\cup\{\Omega_j^*g_{jk}\}_{j\in\tau_0, k\in K_j\setminus \sigma_j}\cup\{\Omega_j^*g_{jk}\}_{j\in\mathbb{N}\setminus\tau_0, k\in K_j}$ is a $\Theta$-frame for $\mathcal{H}$ with lower frame bound $A$ and upper frame bound $B$.
 \end{defn}
A counter-example (Example 3.3 in \cite{XZZ}) by Xiao et al. in the context of characterization of weaving of $\Theta$-frames (see Theorem \ref{thm1}) is incorrect. Further, there is no new characterization of weaving of $g$-frames in \cite{XZZ}. This is given in the following remarks.
\begin{rem}
 Note that in Definition \ref{def3}, $I_j=Q_j=K_j$. Thus, weaving of $\Theta$-frames $\{\Lambda_j^*f_{jk}\}_{j\in\mathbb{N}, k\in I_j}$  and $\{\Omega_j^*g_{jk}\}_{j\in\mathbb{N}, k\in Q_j}$ for $\mathcal{H}$ given in Definition  \ref{def3} of \cite{XZZ} is meaningless  as the indexing sets $I_j$ and $Q_j$ need not be same.
\end{rem}

\begin{rem}
Due to different indexing sets, we cannot use Definition \ref{def3} in Theorem \ref{thm1} for the weaving of $\Theta$-frames $\{\Lambda_j^*f_{jk}\}_{j\in\mathbb{N}, k\in I_j}$ and $\{\Omega_j^*g_{jk}\}_{j\in\mathbb{N}, k\in Q_j}$.
 \end{rem}
\begin{rem}
 In Example 3.3 of \cite{XZZ}, Xiao et al. showed the existence of woven $g$-frames $\{\Lambda_{j}\}_{j=1}^\infty$  and  $\{\Omega_{j}\}_{j=1}^\infty$ such that the frames $\{\Lambda_j^*f_{jk}\}_{j\in\mathbb{N}, k\in K_j}$ and $\{\Omega_j^*g_{jk}\}_{j\in\mathbb{N}, k\in K_j}$ are not woven. However, they used Definition \ref{def3}  for the weaving of frames $\{\Lambda_j^*f_{jk}\}_{j\in\mathbb{N}, k\in K_j}$ and $\{\Omega_j^*g_{jk}\}_{j\in\mathbb{N}, k\in K_j}$, where $\Theta$ is the identity operator. Using Example 3.3. of \cite{XZZ}, Xiao et al.,  claimed that the part (ii) $\implies$ (i) of Theorem \ref{thm1} is incorrect. But, their claim is wrong as  Definition \ref{def1} was used in Theorem \ref{thm1} while Definition \ref{def3} was used in Example 3.3 of \cite{XZZ}. Hence, Example 3.3 \cite{XZZ} does not contradict Theorem \ref{thm1}. Furthermore, Theorem 3.9 of \cite{DVV} is mathematically correct.
 \end{rem}

 \begin{rem}\label{remII}
  Definition \ref{def1} and Definition \ref{def3} are equivalent if $|I_j|=|Q_j|=|K_j|=1$ for all $j\in\mathbb{N}$.
  \end{rem}
  \begin{proof}
Let $I_j=Q_j=K_j$ for $j\in\mathbb{N}$.

First suppose that the $\Theta$-frames  $\{\Lambda_j^*f_{jk}\}_{j\in\mathbb{N}, k\in I_j}$  and $\{\Omega_j^*g_{jk}\}_{j\in\mathbb{N}, k\in Q_j}$  are woven  with universal bounds $A$ and $B$ (in the sense of Definition \ref{def1}).

Let $\tau_0$ be any subset of $\mathbb{N}$, and $\sigma_j\subseteq K_j$ for $j\in\tau_0$.\\
\textbf{Case(i):} Let $\sigma_j$ is non-empty, that is, $\sigma_j=K_j$. Then, $K_j\setminus\sigma_j$ is an empty set and hence we have
\begin{align*}
&\{\Lambda_j^*f_{jk}\}_{j\in\tau_0, k\in \sigma_j}\cup\{\Omega_j^*g_{jk}\}_{j\in\tau_0, k\in K_j\setminus \sigma_j}\cup\{\Omega_j^*g_{jk}\}_{j\in\mathbb{N}\setminus\tau_0, k\in K_j}\\
&=\{\Lambda_j^*f_{jk}\}_{j\in\tau_0, k\in K_j}\cup\{\Omega_j^*g_{jk}\}_{j\in\mathbb{N}\setminus\tau_0, k\in K_j}\\
&=\{\Lambda_j^*f_{jk}\}_{j\in\tau_0, k\in I_j}\cup\{\Omega_j^*g_{jk}\}_{j\in\mathbb{N}\setminus\tau_0, k\in Q_j}
 \end{align*}
 Thus, $\{\Lambda_j^*f_{jk}\}_{j\in\tau_0, k\in \sigma_j}\cup\{\Omega_j^*g_{jk}\}_{j\in\tau_0, k\in K_j\setminus \sigma_j}\cup\{\Omega_j^*g_{jk}\}_{j\in\mathbb{N}\setminus\tau_0, k\in K_j}$ is a $\Theta$-frame for $\mathcal{H}$ with lower frame bound $A$ and upper frame bound $B$.\\
\textbf{Case(ii):} Let  $\sigma_j$ is empty, that is, $K_j\setminus\sigma_j=K_j$. Then, we have
\begin{align*}
&\{\Lambda_j^*f_{jk}\}_{j\in\tau_0, k\in \sigma_j}\cup\{\Omega_j^*g_{jk}\}_{j\in\tau_0, k\in K_j\setminus \sigma_j}\cup\{\Omega_j^*g_{jk}\}_{j\in\mathbb{N}\setminus\tau_0, k\in K_j}\\
&=\{\Omega_j^*g_{jk}\}_{j\in\tau_0, k\in K_j}\cup\{\Omega_j^*g_{jk}\}_{j\in\mathbb{N}\setminus\tau_0, k\in K_j}\\
&=\{\Omega_j^*g_{jk}\}_{j\in\mathbb{N}, k\in K_j}
\end{align*}
Thus, $\{\Lambda_j^*f_{jk}\}_{j\in\tau_0, k\in\sigma_j}\cup\{\Omega_j^*g_{jk}\}_{j\in\tau_0, k\in K_j\setminus \sigma_j}\cup\{\Omega_j^*g_{jk}\}_{j\in\mathbb{N}\setminus\tau_0, k\in K_j}$ is a $\Theta$-frame for $\mathcal{H}$ with lower frame bound $A$ and upper frame bound $B$. Hence, $\Theta$-frames $\{\Lambda_j^*f_{jk}\}_{j\in\mathbb{N}, k\in I_j}$  and $\{\Omega_j^*g_{jk}\}_{j\in\mathbb{N}, k\in Q_j}$ (according to Definition \ref{def3}) are woven with universal bounds $A$ and $B$.

\vspace*{10pt}

In the other direction, suppose that  $\Theta$-frames $\{\Lambda_j^*f_{jk}\}_{j\in\mathbb{N}, k\in I_j}$  and $\{\Omega_j^*g_{jk}\}_{j\in\mathbb{N}, k\in Q_j}$ are woven with universal bounds $A$ and $B$ (in the sense of Definition \ref{def3}).

Let $\sigma$ be any subset of $\mathbb{N}$. Take $\sigma_j=K_j$ for $j\in\sigma$. Then, $K_j\setminus \sigma_j$ is an empty set and hence
\begin{align*}
&\{\Lambda_j^*f_{jk}\}_{j\in\sigma, k\in I_j}\cup\{\Omega_j^*g_{jk}\}_{j\in\mathbb{N}\setminus\sigma, k\in Q_j}\\
&=\{\Lambda_j^*f_{jk}\}_{j\in\sigma, k\in K_j}\cup\{\Omega_j^*g_{jk}\}_{j\in\mathbb{N}\setminus\sigma, k\in K_j}\\
&=\{\Lambda_j^*f_{jk}\}_{j\in\sigma, k\in \sigma_j}\cup\{\Omega_j^*g_{jk}\}_{j\in\sigma, k\in K_j\setminus \sigma_j}\cup\{\Omega_j^*g_{jk}\}_{j\in\mathbb{N}\setminus\sigma, k\in K_j}.
\end{align*}
Thus, $\{\Lambda_j^*f_{jk}\}_{j\in\sigma, k\in I_j}\cup\{\Omega_j^*g_{jk}\}_{j\in\mathbb{N}\setminus\sigma, k\in Q_j}$ is a $\Theta$-frame for $\mathcal{H}$ with upper and lower frame bounds $A$ and $B$, respectively. Hence, $\{\Lambda_j^*f_{jk}\}_{j\in\mathbb{N}, k\in I_j}$  and $\{\Omega_j^*g_{jk}\}_{j\in\mathbb{N}, k\in Q_j}$ are woven with universal bounds $A$ and $B$ (in the sense of Definition \ref{def1}).
\end{proof}

By invoking Remark \ref{remII}, we have the following corollary to the Theorem \ref{thm1}.
\begin{cor}\label{cor1}
Suppose that $\Lambda \equiv\{\Lambda_{j}\}_{j=1}^\infty$  and  $\Omega \equiv\{\Omega_{j}\}_{j=1}^\infty$ are $\Theta$-$g$-frames for $\mathcal{H}$ with respect to $\{ \mathcal{H}_j\}_{j=1}^{\infty}$ and $\{ \mathcal{W}_j\}_{j=1}^{\infty}$, respectively. Let $\{f_{jk}\}_{k\in K_j \subset \mathbb{N}}$ and $\{g_{jk}\}_{k\in K_j \subset \mathbb{N}}$ be orthonormal bases for $\mathcal{H}_j$ and $\mathcal{W}_j$, respectively $(j \in \mathbb{N})$. If $|K_j|=1$ for all $j\in\mathbb{N}$, then the following conditions are equivalent.
\begin{enumerate}[$(i)$]
 \item $\Lambda$ and $\Omega$  are $\Theta$-$g$-woven.
 \item  $\{\Lambda_j^*f_{jk}\}_{j\in\mathbb{N}, k\in K_j}$ and $\{\Omega_j^*g_{jk}\}_{j\in\mathbb{N}, k\in K_j}$ are woven $\Theta$-frames for $\mathcal{H}$ (weaving in the sense of Definition \ref{def3}).
 \end{enumerate}
\end{cor}

\begin{rem}
One can obtain Theorem 3.4 of \cite{XZZ} (in which  Xiao et al. claim new characterization of weaving $g$-frames)  by taking $\Theta$ to be the identity operator in Corollary \ref{cor1}.
\end{rem}

%
%

\end{document}